\documentclass[11pt,oneside]{article}

\usepackage[utf8]{inputenc} 
\usepackage[T1]{fontenc}    
\usepackage[english]{babel}

\usepackage{parskip}
\usepackage{geometry}
\RequirePackage{endnotes}

\usepackage[linesnumbered,ruled,vlined]{algorithm2e}
\usepackage{algpseudocode}
\usepackage{tikz}
\usepackage{array}
\newcolumntype{L}[1]{>{\RaggedRight\arraybackslash}m{#1}} 
\usepackage[sort&compress,numbers]{natbib}

\usepackage{mathtools}
\usepackage{amssymb,amsmath,amsfonts}
\usepackage{xcolor}
\usepackage{bm}

\usepackage{verbatim}
\usepackage[colorlinks=true,citecolor=blue]{hyperref}
\usepackage{url}
\usepackage[capitalize]{cleveref}
\crefformat{equation}{(#2#1#3)}
\crefrangeformat{equation}{(#3#1#4)--(#5#2#6)}
\usepackage{booktabs}

\usepackage{multirow}
\usepackage{booktabs}
\usepackage{pgfplots} 
\usepackage{subcaption} 
\usepgfplotslibrary{groupplots}

\definecolor{gold}{rgb}{0.85,0.65,0}
\colorlet{dgreen}{green!60!black}

\newcommand{\nomm}{\texttt{n-OMM}\xspace}
\newcommand{\eomm}{\texttt{e-OMM}\xspace}
\newcommand{\ceomm}{\texttt{ce-OMM}\xspace}
\newcommand{\romma}{\texttt{ROMMA}\xspace}
\newcommand{\aromma}{\texttt{a-ROMMA}\xspace}
\newcommand{\pumma}{\texttt{PUMMA}\xspace}
\newcommand{\pctr}{\texttt{Perceptron}\xspace}
\newcommand{\alma}{\texttt{ALMA}\xspace}

\newcommand{\grad}{\ensuremath{\nabla}}

\let\emptyset\varnothing
\newcommand{\set}[1]{\left\{#1\right\}}

\newcommand{\barD}{\ensuremath{\overline{D}}}

\def\bbN{{\mathbb{N}}}

\def\bbR{{\mathbb{R}}}

\def\cA{\mathcal{A}}

\DeclarePairedDelimiterX{\inner}[2]{\langle}{\rangle}{#1, #2}

\DeclareMathOperator{\conv}{conv}

\DeclareMathOperator{\sign}{sign}

\usepackage{amsthm}

\newcommand{\epr}{\hfill$\blacksquare$}
\newtheorem{theorem}{Theorem}[section]
\newtheorem{lemma}{Lemma}[section]

\newtheorem{assumption}{Assumption}
\crefname{assumption}{assumption}{assumptions}
\Crefname{assumption}{Assumption}{Assumptions}
\theoremstyle{definition}
\newtheorem{remark}{Remark}[section]
\newtheorem{example}{Example}[section]
\newtheorem{definition}{Definition}[section]

\def\proof#1{{\it #1\enskip }\ignorespaces}
\def\Halmos{\unskip\hfill$\square$}

\DeclareMathOperator*{\argmax}{arg\,max}

\usepackage{authblk}

\begin{document}

\title{Efficient Online Large-Margin Classification via Dual Certificates}

\author[1]{Nam Ho-Nguyen}
\author[2]{Fatma K{\i}l{\i}n\c{c}-Karzan}
\author[1]{Ellie Nguyen}
\author[2]{Lingqing Shen}
\affil[1]{Discipline of Business Analytics, The University of Sydney}
\affil[2]{Tepper School of Business, Carnegie Mellon University}
\date{}

\maketitle

\begin{abstract}
Online classification is a central problem in optimization, statistical learning and data science. Classical algorithms such as the perceptron offer efficient updates and finite mistake guarantees on linearly separable data, but they do not exploit the underlying geometric structure of the classification problem.
We study the offline maximum margin problem through its dual formulation and use the resulting geometric insights to design a principled and efficient algorithm for the online setting. A key feature of our method is its \emph{translation invariance}, inherited from the offline formulation, which plays a central role in its performance analysis. Our theoretical analysis yields improved mistake and margin bounds that depend only on translation-invariant quantities, offering stronger guarantees than existing algorithms under the same assumptions in favorable settings. In particular, we identify a parameter regime where our algorithm makes at most two mistakes per sequence, whereas the perceptron can be forced to make arbitrarily many mistakes. Our numerical study on real data further demonstrates that our method matches the computational efficiency of existing online algorithms, while significantly outperforming them in accuracy.  
\end{abstract}

\section{Introduction}\label{sec:Intro}

Online learning is a topic with rich connections to several fields of interest and many real-world applications such as portfolio selection \citep{LiHoi2014}, personalized recommendation systems \citep{LiEtAl2010}, and online advertising \citep{McMahanEtAl2013}; see \citet{Cesa-BianchiLugosi2006book,HoiEtAl2021} for a more comprehensive overview. Within this domain, online classification
plays a central role: given a sequence of feature vectors, the learner must issue immediate predictions and incrementally update its model as feedback is revealed. This sequential nature makes the setting particularly relevant for modern environments characterized by large-scale, high-dimensional, and high-velocity data streams.

Online classification methods often rely on two structural assumptions: \emph{linearity}, where the decision boundary at each iteration is a hyperplane, and \emph{separability}, where positive and negative points are separated by a margin.
These assumptions not only simplify the problem but also yield valuable geometric structure that can be exploited to design efficient algorithms with provable guarantees.
Kernel methods further enable nonlinear decision boundaries while preserving many of the appealing geometric features. 

The perceptron algorithm, introduced by \citet{rosenblatt_perceptron_1958}, is the first and most widely known online classification method. In the linearly separable case, \citet{novikoff_convergence_1962} proved a finite mistake bound for the perceptron inversely proportional to the margin; the non-separable case was later analyzed by \citet{FreundSchapire1999}.
Many perceptron variants retain the additive update rule but modify the step sizes to improve performance and broaden applicability
\citep{shalev_ballseptron_2005,bshouty_cosine_2016,crammer_mira_2003,CrammerEtAl2006,KivinenEtAl2002}. Kernelized versions enable nonlinear decision boundaries \citep{aizerman_kernelpctr_1964}, which also inspired further enhancements \citep{dekel_forgetron_2005,dekel_forgetron_2008,he_fuseptron_2012,cesa-bianchi_sampling_2006,hoi_omkc_2013,KivinenEtAl2004}.
Despite these advances,
the perceptron and its variants suffer from a central limitation: they do not necessarily recover a large-margin classifier, which is essential for robustness to perturbations \citep[see, e.g.,][]{JavanmardMehrabi2024}, improved statistical generalization \citep[Ch.~9]{AnthonyBartlett1999book}, and potentially higher accuracy.

Margin maximization has been studied extensively as well, first in the offline setting, starting with linear and quadratic programming approaches by \citet{Mangasarian1965} and \citet{Rosen1965}.
\citet{BoserEtAl1992} later recognized that the quadratic programming formulation is indeed equivalent to the \emph{maximum-margin problem} underlying support vector machines.
\citet{AnlaufBiehl1989} suggested an iterative algorithm for this quadratic program,
which was later extended to the kernel setting by \citet{FriessEtAl1998}. Around the same time, \citet{Kowalczyk2000} and \citet{KeerthiEtAl2000} proposed other iterative algorithms with extensions to both kernel and non-separable cases.
The algorithm of \citet{Kowalczyk2000} can be adapted to the online setting, but obtaining guarantees requires knowledge of problem parameters, such as margin size, that are not available a priori.
The algorithms ALMA \citep{gentile_new_2000}, ROMMA \citep{li_relaxed_1999} and PUMMA \citep{IshibashiEtAl2008} provide implementable margin guarantees in the online setting. ALMA achieves any chosen $\rho$-fraction of a version of the margin where data are normalized; ROMMA matches ALMA’s guarantee without requiring $\rho$ and can converge to the perfect margin $\gamma_*$; and PUMMA extends ROMMA to arbitrary $p$-norms for $p \geq 2$, as well as nonhomogeneous hyperplanes. 
See \cref{tab:perf_guarantees} for a comparison of the performance guarantees of these algorithms; in particular, their existing mistake bound guarantees match the perceptron.

It is well-known that the perceptron's mistake bound, shared by many of these algorithms, is tight \citep[see, e.g.,][Exercise 8.3]{mohri_foundations_2012}. This classical tight example shows that as the problem dimension $d$ grows, the margin decreases, and \emph{any} deterministic algorithm is forced to make $d$ mistakes. However, such guarantees capture only worst-case behavior and are too coarse to differentiate the performance of algorithms in more typical scenarios. In particular,
all existing performance guarantees for these algorithms are expressed in \emph{translation-dependent} terms, and the algorithms themselves (with the exception of PUMMA) are also translation dependent.
Consequently, prior analyses have not fully leveraged the geometric structure of maximum-margin classifiers, nor have they produced geometry-based mistake bounds that can surpass the classical barrier in favorable settings.

In this paper, we identify \emph{translation invariance} (formally defined in \cref{rem:translation-invariance}), a fundamental property of the offline maximum-margin problem, as a key ingredient that can be leveraged to design online classification algorithms achieving stronger guarantees in favorable settings, while matching worst-case guarantees up to a logarithmic factor.
By carefully analyzing the dual formulation of the offline maximum margin problem, we propose a new efficient 
\emph{translation-invariant} algorithm for online classification in the linearly separable case. 
The crux of our analysis is based on a geometric interpretation of the \emph{dual certificates}, and a bound on how much the margin changes as new points are introduced. We use these certificates to identify a principled way to approximate a na\"{i}ve and computationally impractical online maximum margin method, and as a result derive a computationally efficient version. 
Moreover, our approach retains the underlying geometric structure and yields performance guarantees
on both the number of mistakes and different scales of margin violations expressed solely on translation-invariant quantities.
Through theoretical comparisons, illustrative examples, and numerical results, we demonstrate the benefit of exploiting the translation-invariant structure over existing methods.
A central outcome of our refined analysis {is the identification,
for the first time in the literature, of a parameter regime where there exists a large theoretical gap in mistake bounds between the translation-dependent algorithms and ours. This finding thus shows that the classical mistake bound guarantees of the existing algorithms, such as perceptron, can be substantially improved in favorable cases. Specifically, in parameter regimes where the positive points and negative points are well clustered and these two clusters are well separated, we show that
our algorithm makes at most two mistakes, whereas an algorithm such as the perceptron makes a large number (see \cref{sec:comparison} and \cref{ex:bad-perceptron}). On the practical side, our algorithm requires only the choice of a norm, from a very broad class, and an optional, easily interpretable hyper-parameter.
Collectively, our results contribute new geometric, optimization-based insights and practical improvements for the online classification problem.

We formally describe the online classification framework, along with our notation, in \cref{sec:problem}. In \cref{sec:online-max-margin}, we present our algorithm, which constructs online classifiers by approximating the offline maximum margin solution using dual certificates. Our performance analysis in \cref{sec:performance-guarantees} establishes guarantees on both the margin and the number of mistakes based on translation-invariant quantities, reflecting the geometric structure inherited from the offline problem. In \cref{sec:implementation}, we compare our guarantees with those of existing methods and provide a concrete example illustrating the practical impact of our approach. Finally, \cref{sec:numerical} presents an empirical study on real-world data, highlighting the effectiveness and robustness of our algorithms (and thus translation invariance) in large-scale settings.

\section{Problem description}\label{sec:problem}

Let $\|\cdot\|$ be a given norm on $\bbR^d$, and let $\|\cdot\|_*$ denote its dual norm.
We will work with the following assumption on the norm.
\begin{assumption}\label{ass:norm}
The norm $\|\cdot\|$ is strictly convex and differentiable away from $0$. Consequently, the dual norm $\|\cdot\|_*$ is also strictly convex and differentiable away from $0$.
\end{assumption}
For $x,w \neq 0$, denote $\ell(x) := \grad \|x\|$ and $\ell_*(w) := \grad \|w\|_*$ to be the gradient of the primal and dual norms, respectively. From convex analysis, we have $\|\ell(x)\|_* = \|\ell_*(w)\| = 1$, and $\ell(x)^\top x = \|x\|$, $\ell_*(w)^\top w = \|w\|_*$. Furthermore, when $\|x\| > 0$ and $\|w\|_* > 0$, we have
$\ell_*(\ell(x)) = x/\|x\|$ and $\ell(\ell_*(w)) = w/\|w\|_*$. 

We first describe the \emph{online linear classification} problem. At each time $t=1,2,\ldots$, the learner receives a feature vector $x_t \in \bbR^d$, and predicts its binary label using a linear classifier $(w_t,b_t) \in \bbR^d \times \bbR$ via $\hat{y}_t := \sign(w_t^\top x_t + b_t) \in \{\pm 1\}$. Throughout the paper, we will assume that $\sign(0)=1$.
The true label $y_t \in \{\pm 1\}$ is then revealed. Based on the feedback, the learner then updates the classifier $(w_{t+1},b_{t+1})$ for the next time step.

For a given data pair $(x,y)$ and a classifier $(w,b)$, the \emph{margin} is defined to be the smallest distance by which we need to perturb the feature vector $x$ in order to cause a misclassification by the classifier $(w,b)$, which admits the following closed-form expression:
\[ \gamma(x,y; w,b) = \frac{\max\set{ 0,\, y(w^\top x + b) }}{\|w\|_*}. \]
Throughout, we will make the following linear separability and boundedness assumptions on the data.
\begin{assumption}\label{ass:margin-diameter}
There exists an (unknown) \emph{perfect classifier} $(w_*,b_*)$ such that $y_t = \sign(w_*^\top x_t + b_*)$ for all $t$, and among all such correct classifiers, $(w_*,b_*)$ achieves the largest margin $\gamma_*$. On this classifier, the margin of each feature vector is positive:
\[ \gamma_* := \inf_t \gamma(x_t,y_t; w_*,b_*) > 0. \]
Furthermore, the following quantities are finite:
\begin{subequations}\label{eq:diameter}
\begin{align}
D_+ &:= \sup_{t,t'}\set{ \|x_t-x_{t'}\| :~ y_t = y_{t'} = +1 },\\
D_- &:= \sup_{t,t'}\set{ \|x_t-x_{t'}\| :~ y_t = y_{t'} = -1 },\\
D &:= \max\{ D_+,\, D_-\}, ~~  \barD := \sup_t \|x_t\|, ~~ r := D/\gamma_*.
\end{align}
\end{subequations}
\end{assumption}
Here, $D_+$ and $D_-$ measure the diameter of the clusters of positive and negative data points, respectively, while $\barD$ is the maximum length of all (both positive and negative) points. Finally, the parameter $r$ is a measure of how well separated the two clusters of positive and negative points are with respect to the margin.

In online classification, one of the most fundamental performance measures of an algorithm is how its number of mistakes grows as a function of $t$. As another performance measure, we will also count the number of times the predicted margin is less than a factor $\phi < 1$ of $\gamma_*$. 
We will make these notions precise in \cref{sec:performance-guarantees} (cf.\ \cref{eq:m-phi-counts}).

\begin{remark}[Translation invariance]\label{rem:translation-invariance}
Consider two online classification problems, where one receives the sequence of points is $\{x_t,y_t\}_{t \geq 1}$ and the other receives \emph{translated} feature vectors $\{x_t+u,y_t\}_{t \geq 1}$, where $u$ is some fixed vector. The parameters $D_+,D_-,D,r$ are identical in both problems, and if $(w_*,b_*)$ is a perfect classifier for the first problem with margin $\gamma_*$, then $(w_*,b_*-w_*^\top u)$ is a perfect classifier for the second problem with the same margin $\gamma_*$. In this sense, the two problems are essentially the same, thus we expect a good algorithm for the online classification problem to be invariant to translation by $u$. 

While the parameters $D$, $r$ and $\gamma_*$ are all translation invariant, the parameter $\barD$ is translation dependent.
Existing online classification algorithms admit performance guarantees that are given in terms of $\gamma_*$ and $\barD$ (see \cref{tab:perf_guarantees}), and hence may inherit the translation-dependent property of the latter parameter. In contrast, $D$, $r$ and $\gamma_*$ are all translation invariant, so any algorithm admitting guarantees depending solely on these parameters
will naturally be translation-invariant as well.
As we will see, our \cref{alg:max-margin} defined below indeed has this property (see \cref{rem:translation-invariance-alg}), and our analysis depends only on $D$, $r$ and $\gamma_*$ (see \cref{thm:margin-performance,sec:implementation-l2}).
This is in sharp contrast to the translation-dependent guarantees of the previous algorithms.
Moreover, we will show that the translation dependence of an algorithm can have a detrimental effect on its performance analytically in \cref{ex:bad-perceptron} and computationally in \cref{sec:numerical}.
\epr
\end{remark}

\section{The online maximum margin algorithm}\label{sec:online-max-margin}

Before describing our algorithm, we provide some fundamental results on maximum margin classifiers. Consider two finite linearly separable sets $V_+,V_- \subset \bbR^d$, and a labeling function $y(\cdot)$ such that $y(v)=+1$ for $v \in V_+$ and $y(v) = -1$ for $v \in V_-$. The \emph{(offline) maximum margin problem} aims to find a linear separator $(w,b)$ by solving
\begin{equation}\label{eq:margin}
\gamma(V_+,V_-) := \max_{w,b} \min_{v \in V_+ \cup V_-} \gamma(v,y(v); w,b),
\end{equation}
i.e., by maximizing the distance from $v$ to the misclassification hyperplane $H_v := \left\{ s : y(v)(w^\top s + b) \leq 0 \right\}$. 

Our design and analysis of the online maximum margin algorithm is based on the connection between \eqref{eq:margin} and the following convex optimization problem:
\begin{equation}\label{eq:convex-hull-distance}
\tau(V_+,V_-) := \min_{v_+,v_-} \set{ \|v_+ - v_-\| :~ \begin{aligned} &v_+ \in \conv(V_+),\\ &v_- \in \conv(V_-) \end{aligned} }.
\end{equation}
\begin{definition}[Dual certificates]\label{def:dual-certificate}
Given finite linearly separable sets $V_+,V_- \subset \bbR^d$, let $v_+ \in \conv(V_+), v_- \in \conv(V_-)$ be solutions to \cref{eq:convex-hull-distance}. We call $v_+,v_-$ the \emph{dual certificates} of \cref{eq:margin} as they certify the distance between $\conv(V_+)$ and $\conv(V_-)$. We sometimes use the terminology that $v_+,v_-$ are the dual certificates of $\gamma(V_+,V_-)$.
\epr
\end{definition}

\begin{remark}\label{rem:translation-invariance-offline}
Both the offline maximum margin problem \cref{eq:margin} and the convex problem \cref{eq:convex-hull-distance} are translation invariant: if the positive and negative sets are $V_+ + \{u\}$ and $V_- + \{u\}$, then the solution of \cref{eq:margin} is simply $(w,b-w^\top u)$ with the same optimal margin $\gamma(V_+ + \{u\}, V_- + \{u\}) = \gamma(V_+,V_-)$ and the dual certificates $v_+ + u$ and $v_- + u$.
\epr
\end{remark}

The dual certificate terminology stems from the fact that \cref{eq:convex-hull-distance} is derived by examining the dual of \cref{eq:margin}, and that the optimal solution of \cref{eq:margin} can be directly obtained from the certificates. We make this precise in the following result.
\begin{lemma}\label{lemma:separable-max-margin}
Suppose $V_+,V_- \subset \bbR^d$ are finite linearly separable sets. Let $v_+ \in \conv(V_+)$ and $v_- \in \conv(V_-)$ be the dual certificates, so that $\|v_+ - v_-\| = \tau(V_+,V_-)$. Then, $\gamma(V_+,V_-) = \tau(V_+,V_-)/2$, and an optimal solution $(w,b)$ for \cref{eq:margin} is given by 
$ w := \ell\left( v_+ - v_- \right)$ and $b := -\frac{1}{2} w^\top (v_+ + v_-).
$ \end{lemma}
\begin{remark}\label{rem:hard-margin-SVM}
\citet{BennettBredensteiner2000} showed the equivalence of \cref{eq:margin} and \cref{eq:convex-hull-distance} for $\ell_2$-norm; our \cref{lemma:separable-max-margin} essentially generalizes this to \emph{any} norm satisfying \cref{ass:norm}. Beyond this generalization, the real value of \cref{lemma:separable-max-margin} lies in the geometric interpretation it supplies: the dual certificates $v_+,v_-$ act as compact representatives of the positive and negative classes, capturing exactly the information needed to characterize the maximum margin. This insight gives us an important leverage in algorithm design as follows. Instead of storing and processing the entire stream of past examples, \cref{alg:max-margin} maintains only the current representatives $v_+,v_-$, updating them at every round. This design thus allows us to keep per-iteration complexity low while capturing the geometric essence of \cref{eq:margin}, which is key to our performance analysis in \cref{sec:performance-guarantees}.
\epr
\end{remark}
\proof{Proof of \cref{lemma:separable-max-margin}.}
By expanding problem \cref{eq:margin} we arrive at
\begin{align*}
\max_{w,b,\alpha} \set{ \alpha: ~
\begin{array}{l}
 \frac{\max\{0, w^\top v + b\}}{\|w\|_*} \geq \alpha, \quad \forall v \in V_+, \\
\frac{\max\{0, -(w^\top v + b)\}}{\|w\|_*} \geq \alpha, \quad \forall v \in V_-
\end{array} }.
\end{align*}
Since $V_+,V_-$ are linearly separable, the optimal value of \cref{eq:margin} is $> 0$. Therefore, we can transform the variables $w/(\alpha \|w\|_*) \to w$ and $b/(\alpha \|w\|_*) \to b$ to arrive at the well-known exact convex reformulation of \cref{eq:margin}, i.e., the so-called \emph{hard-margin support vector machine} problem:
\begin{align*}
\min_{w,b} \set{\|w\|_*:~\begin{array}{l}
 w^\top v + b \geq 1, \quad \forall v \in V_+\\
 w^\top v + b \leq -1, \quad \forall v \in V_-
\end{array}} .
\end{align*}
Furthermore, $\|w\|_* = 1/\alpha$ in the new reformulation, thus if the solution $(w,b)$ solves the convex reformulation, then $\gamma(V_+,V_-) = 1/\|w\|_*$.

Let us now examine the conic dual of the preceding reformulation given by 
\begin{align*}
\max_{\lambda,\eta} \quad & \sum_{v \in V_+} \lambda_v + \sum_{v \in V_-} \eta_v\\
\text{s.t.} \quad & \lambda, \eta \geq 0\\
& \left\| \sum_{v \in V_+} \lambda_v v - \sum_{v \in V_-} \eta_v v \right\| \leq 1\\
& \sum_{v \in V_+} \lambda_v = \sum_{v \in V_-} \eta_v.
\end{align*}
Now, given some solution $\lambda,\eta$ to this conic dual, we define $c := \sum_{v \in V_+} \lambda_v = \sum_{v \in V_-} \eta_v$ and our dual certificates $v_+,v_-$ as $v_+ := {1\over c} (\sum_{v \in V_+} \lambda_v v) \in \conv(V_+)$ and $v_- := {1\over c} (\sum_{v \in V_-} \eta_v v) \in \conv(V_-)$.
Then, we have $\left\| \sum_{v \in V_+} \lambda_v v - \sum_{v \in V_-} \eta_v v \right\| = c\|v_+ - v_-\| \leq 1$.
However, notice that $c$ is essentially independent of $v_+,v_-$, except for the constraint that $c \leq 1/\|v_+ - v_-\|$, therefore the dual can be written as
\begin{align*}
\max_{c,v_+,v_-} \set{
2c:~ \begin{array}{l}
c \leq 1/\|v_+ - v_-\|, \\
v_+ \in \conv(V_+), \ v_- \in \conv(V_-)
\end{array}
} . 
\end{align*}
The objective function is simply $2c$, so to maximize it we set $c = 1/\|v_+ - v_-\|$, and choose $v_+,v_-$ to minimize $\|v_+ - v_-\|$, which is exactly \cref{eq:convex-hull-distance}. This implies that if $(w,b)$ is optimal for the convex reformulation, then $\|w\|_* = 2/\tau(V_+,V_-)$. Since $1/\|w\|_* = \gamma(V_+,V_-)$, we have $\gamma(V_+,V_-) = \tau(V_+,V_-)/2$.

Considering again the convex reformulation and its dual, the KKT optimality conditions state that
\begin{align*}
w^\top \left( \sum_{v \in V_+} \lambda_v v - \sum_{v \in V_-} \eta_v v \right) &= \|w\|_* , \\
\left\| \sum_{v \in V_+} \lambda_v v - \sum_{v \in V_-} \eta_v v \right\| &= 1 , \\
\lambda_v(1-(w^\top v +b)) &= 0, \ \forall v \in V_+ ,\\
\eta_v(-1 - (w^\top v +b)) &= 0, \ \forall v \in V_-.
\end{align*}
The first two of these conditions together imply that $\sum_{v \in V_+} \lambda_v v - \sum_{v \in V_-} \eta_v v = \argmax_{x : \|x\| \leq 1} w^\top x = \ell_*(w)$.
Using the relationship between $\ell$ and $\ell_*$, we have $\ell(\ell_*(w)) = \frac{w}{\|w\|_*} = \ell\left( \sum_{v \in V_+} \lambda_v v - \sum_{v \in V_-} \eta_v v \right)
= \ell(v_+ - v_-)$.
Summing the third condition over $v \in V_+$ implies that 
$\sum_{v \in V_+} \lambda_v(1-w^\top v) = b\sum_{v \in V_+} \lambda_v. $ Then, dividing by $\sum_{v \in V_+} \lambda_v$ we get $b = 1 - w^\top v_+$. The same argument for the fourth condition gives $b = -(1 + w^\top v_-)$. Taking the average gives $b = -\frac{1}{2} w^\top (v_+ + v_-)$. Then $(w,b)$ is optimal for the convex reformulation, and thus $(w/\|w\|_*, b/\|w\|_*)$ is also optimal for \cref{eq:margin} since it is scale-invariant.
\Halmos
\endproof

\cref{alg:max-margin} formally describes two variants of our proposed method: a \emph{na\"{i}ve} implementation that at each iteration simply keeps all points and solves the maximum margin problem with all the data seen thus far; and an \emph{efficient} implementation that utilizes the idea discussed in \cref{rem:hard-margin-SVM} of using the dual certificates $v_+,v_-$ as representative points. By exploiting the fact that the dual certificates adequately capture the geometry of the maximum margin problem, in \cref{sec:performance-guarantees} we in fact show  that both versions enjoy the \emph{same} guarantees.
\begin{algorithm}[tb]
\caption{Online Maximum Margin}
\label{alg:max-margin}
\KwIn{Norm $\|\cdot\|$ along with its dual $\|\cdot\|_*$, aggressiveness parameter $\rho \in [0,1]$.}
\KwOut{Sequence of linear classifiers $\{(w_t, b_t)\}_{t \geq 1}$.}

Obtain $x_1$, predict $\hat{y}_1 = +1$, receive $y_1$\;

Initialize $V_{3,+} := \{x_1\},V_{3,-} := \emptyset$ if $y_1 = +1$; otherwise initialize $V_{3,+} := \emptyset, V_{3,-} := \{x_1\}$ if $y_1 = -1$.

\While{one of $V_{3,+},V_{3.-}$ is empty}{
    Obtain $x_t$, predict $\hat{y}_t = y_1$, receive $y_t$\;
    \If{$y_t \neq \hat{y}_t$}{
        Add $x_t$ to $V_{3,+}$ if $y_t = +1$; otherwise
        $y_t=-1$ so add $x_t$ to $V_{3,-}$\;
    }
}

Reset iteration counter to $t=2$. Since by construction $V_{3,+}$ and $V_{3,-}$ are singletons, denote $v_{3,+},~v_{3,-}$ as the elements. Then $v_{3,+},v_{3,-}$ are trivially the dual certificates of $\gamma_3 := \gamma(V_{3,+},V_{3,-})$\;
Set $w_3 := \ell(v_{3,+}-v_{3,-})$ and $b_3 := -\frac{1}{2} w_3^\top (v_{3,+} + v_{3,-})$\;

\For{$t = 3, 4, \ldots$}{
    Obtain $x_t$, predict $\hat{y}_t \leftarrow \sign(w_t^\top x_t + b_t)$, receive $y_t$, and compute $a_t := y_t (w_t^\top x_t + b_t)$\;

    \If{$a_t < \rho \gamma_t$ (update occurs)}{
    \If{Na\"{i}ve implementation}{
        Update $V_{t+1,+} = V_{t,+} \cup \{x_t\}$ and $V_{t+1,-} = V_{t,-}$ if $y_t = 1$; otherwise
        $y_t = -1$ so set $V_{t+1,+} = V_{t,+}$ and $V_{t+1,-} = V_{t,-} \cup \{x_t\}$\;
    }\ElseIf{Efficient implementation}{
        Update $V_{t+1,+} = \{v_{t,+}, x_t\}$ and $V_{t+1,-} = \{v_{t,-}\}$ if $y_t = +1$; otherwise
        $y_t = -1$ so set $V_{t+1,+} = \{v_{t,+}\}$ and $V_{t+1,-} = \{v_{t,-}, x_t\}$\;
    }
    Solve \eqref{eq:margin} (or \cref{eq:convex-hull-distance}) with $V_{t+1,+},V_{t+1,-}$, and denote the dual certificates as $v_{t+1,+},v_{t+1,-}$\;
        Set $\gamma_{t+1} := \gamma(V_{t+1,+},V_{t+1,-})$, $w_{t+1} := \ell(v_{t+1,+} - v_{t+1,-})$ and  $b_{t+1} := -\frac{1}{2} w_{t+1}^\top (v_{t+1,+} + v_{t+1,-})$\;
    }\Else{
        Set $v_{t+1,+} := v_{t,+}$, $V_{t+1,+} := V_{t,+}$, $v_{t+1,-} := v_{t,-}$, $V_{t+1,-} := V_{t,-}$, $w_{t+1} := w_t$, $b_{t+1} := b_t$, and $\gamma_{t+1} := \gamma_t$\;
    }
    
    }
\end{algorithm}

\begin{remark}[Initialization and aggressiveness parameter]
The initialization phase of \cref{alg:max-margin} (before the ``for'' loop) is to guarantee that \cref{eq:convex-hull-distance} is well-defined by ensuring that there is at least one positive and one negative point in the sets $V_{3,+}$ and $V_{3,-}$, respectively.

The aggressiveness parameter $\rho \in [0,1]$  controls how often \cref{eq:margin} is solved. Setting $\rho = 0$ means that updates are performed only when mistakes occur, i.e., $a_t<0$, so the algorithm is \emph{conservative}. Setting $\rho = 1$ means that the update occurs as often as possible, and is equivalent to removing the condition $a_t < \rho\gamma_t$ completely, since $(w_t,b_t)$ is already optimal for the corresponding problem \cref{eq:margin} when $a_t \geq \gamma_t$. 
\epr
\end{remark}

\begin{remark}[Translation invariance of \cref{alg:max-margin}]\label{rem:translation-invariance-alg}
Elaborating on \cref{rem:translation-invariance}, it is easy to see that \cref{alg:max-margin} is translation invariant, since it solves \cref{eq:convex-hull-distance} at each iteration which itself is translation invariant by \cref{rem:translation-invariance-offline}. More precisely,
suppose \cref{alg:max-margin} returns the sequence of classifiers $(w_t,b_t)$ for a given online classification problem with the sequence of points $(x_t,y_t)$. Then \cref{alg:max-margin} will return the sequence of classifiers $(w_t,b_t - w_t^\top u)$ when solving the translated version of this online problem with the sequence of points $(x_t+u,y_t)$.
\epr
\end{remark}

\begin{remark}[Complexity of \cref{alg:max-margin}]\label{rem:complexity}
The na\"{i}ve implementation involves solving an offline maximum margin problem \cref{eq:margin} with $t$ points at each iteration $t \geq 3$. This admits a convex reformulation; see the proof of \cref{lemma:separable-max-margin}. When the norm $\|\cdot\|$ is an $\ell_p$-norm, it can be solved via off-the-shelf conic solvers. However, the number of constraints in this problem increases with $t$, and thus this is not efficient. 
In contrast, our new efficient implementation requires solving \cref{eq:margin} with only \emph{three} points at each iteration. Off-the-shelf conic solvers can be used to solve this, but since there are only three points, a closed-form solution can be derived in certain cases, e.g., when $\ell_2$-norm is used. We will describe this further in \cref{sec:implementation}.
\epr
\end{remark}

In the next section, we derive performance guarantees for both versions of \cref{alg:max-margin}.

\section{Performance guarantees}\label{sec:performance-guarantees}

We will show that, as \cref{alg:max-margin} progresses, a certain performance metric associated with it is bounded. To define our performance metric, we first establish a basic fact on the margin $\gamma_t$ of the classifier produced by \cref{alg:max-margin} at iteration $t$.
\begin{lemma}\label{lemma:gamma_bound}
For $t \geq 3$, let $\gamma_t$ be defined as in \cref{alg:max-margin}. Then we have $\gamma_t \geq \gamma_* > 0$.
\end{lemma}
\proof{Proof of \cref{lemma:gamma_bound}.}
Recall that by definition, the optimal classifier $(w_*,b_*)$  satisfies $y_k(w_*^\top x_k + b_*)/\|w_*\|_* \geq \gamma_*$ for all $k \geq 1$. Each update iteration (i.e., $a_t < \rho \gamma_t)$ of the na\"{i}ve implementation finds the largest $\gamma$ for which there exists some $(w,b)$ such that $y_k(w^\top x_k + b)/\|w\|_* \geq \gamma$ for
all $k \leq t$ such that $a_k < \rho \gamma_k$.
This means that $\gamma_{t+1} \geq \gamma_*$. For the efficient implementation, in each update iteration, we have that $v_{t,+},v_{t,-}$ are convex combinations of the points $x_k$ for $k \leq t$. Therefore, $(w_*^\top v_{t,+} + b_*)/\|w_*\| \geq \gamma_*$ and $-(w_*^\top v_{t,-} + b_*)/\|w_*\| \geq \gamma_*$. Also, $y_t(w_*^\top x_t + b_*)/\|w_*\| \geq \gamma_*$, so when we solve the maximum margin algorithm on the three points $(v_{t,+},+1),(v_{t,-},-1),(x_t,y_t)$, we get $\gamma_{t+1} \geq \gamma_*$.
\Halmos
\endproof

Now, given $\phi < 1$ we count two types of ``violations'' over the iterations as follows:
\begin{subequations}\label{eq:m-phi-counts}
\begin{align}
\bar{m}(\phi) &:= \#\set{ t \geq 3 :~ \gamma(x_t,y_t;w_t,b_t) \leq \phi \gamma_* },\\
m(\phi) &:= \#\set{ t \geq 3 :~ \gamma(x_t,y_t;w_t,b_t) \leq \phi \gamma_t },
\end{align}
where given a set $\cA$, we define $\#\cA$ to be the cardinality of the set.
\end{subequations}
In words, $\bar{m}(\phi)$ counts the number of times the algorithm encounters a point for which the margin $\gamma(x_t,y_t;w_t,b_t)$ of the current classifier $(w_t,b_t)$ is less than a factor of $\phi$ of the optimal margin $\gamma_*$. In contrast, $m(\phi)$ counts the number of times that the margin $\gamma(x_t,y_t;w_t,b_t)$ is less than $\phi \gamma_t$, where $\gamma_t$ is the predicted margin at time $t$. Thus, while $\bar{m}(\phi)$ depends on $\gamma_*$ hence cannot be computed in practice, $m(\phi)$ can be computed as we run \cref{alg:max-margin}. \cref{lemma:gamma_bound} states that $\gamma_t \geq \gamma_*$ for \cref{alg:max-margin}, hence $\bar{m}(\phi) \leq m(\phi)$. Furthermore, our analysis naturally gives rise to a bound on $m(\phi)$ rather than $\bar{m}(\phi)$. Note that $m(0) = \bar{m}(0)$ is simply the number of mistakes our algorithm makes. There is a small subtlety in \cref{eq:m-phi-counts} in that we count the violations for $t \geq 3$ only. This is due to the initialization phase of \cref{alg:max-margin}, before the ``for'' loop. During the initialization phase, at most two mistakes are made, so the total number of mistakes throughout the algorithm is at most $m(0) + 2$.

Our analysis proceeds as follows. \cref{lemma:max-margin-increment,lemma:kappa} will show that at steps where $a_t \leq \phi \gamma_t$ occurs, the next margin $\gamma_{t+1}$ will decrease by a fixed factor compared to $\gamma_t$, where the factor is defined as the $\kappa_{\circ}$ function in \cref{eq:kappa0}. However, by \cref{lemma:gamma_bound} all margins $\gamma_t$ are bounded below by $\gamma_*$, hence such decreases can occur only finitely many times.

The $\kappa_\circ$ function is defined as follows:
\begin{align}\label{eq:kappa0}
\kappa_\circ\left( \delta, \eta \right) := \max_{u,z} \set{ \min_{\beta \in [0,1]} \left\| u - \beta z \right\| : \begin{aligned} \|u\| &=1\\
\|z\| &\leq \delta\\ \ell(u)^\top z &\geq \eta \end{aligned} }.
\end{align}
\cref{lemma:max-margin-increment} states that when a single point is added to an offline maximum margin problem $\gamma(V_+,V_-)$, the change in the margin can be bounded by the function $\kappa_\circ$ as well as the dual certificates. 
\begin{lemma}\label{lemma:max-margin-increment}
Let $V_+,V_-$ be two finite linearly separable sets, so $\gamma(V_+,V_-) > 0$. Let $v_+ \in \conv(V_+)$, $v_- \in \conv(V_-)$ be the dual certificates of \cref{eq:margin}, and denote $\bar{v} := v_+ - v_-$. Suppose $x$ is a new point such that $V_+ \cup \{x\}$ and $V_-$ are linearly separable. Then, 
\begin{align*}
\frac{\gamma(V_+ \cup \{x\}, V_-)}{\gamma(V_+,V_-)}
&\leq \kappa_\circ\left( \frac{\|v_+ - x\|}{\|\bar{v}\|}, \frac{\ell(\bar{v})^\top(v_+ - x)}{\|\bar{v}\|} \right).
\end{align*}
If $x$ is such that $V_+$ and $V_- \cup \{x\}$ are linearly separable instead, then
\begin{align*}
\frac{\gamma(V_+, V_- \cup \{x\})}{\gamma(V_+,V_-)}
&\leq \kappa_\circ\left( \frac{\|x - v_-\|}{\|\bar{v}\|}, \frac{\ell(\bar{v})^\top(x - v_-)}{\|\bar{v}\|} \right).
\end{align*}
\end{lemma}
\proof{Proof of \cref{lemma:max-margin-increment}.}
Denote $\gamma := \gamma(V_+,V_-)$. Recall from \cref{lemma:separable-max-margin} that $\|v_+ - v_-\| = 2\gamma$. Let $(w,b)$ solve \cref{eq:margin} on $V_+ \cup \{x\}$ and $V_-$, normalized so that $\|w\|_* = 1$, and denote $\bar{\gamma} := \gamma(V_+ \cup \{x\}, V_-)$.
Observe that $w^\top v_+ + b \geq \bar{\gamma}$, $w^\top v_- + b \leq -\bar{\gamma}$, and $w^\top x + b \geq \bar{\gamma}$.
Then, for any $\beta \in [0,1]$ by combining the first and third inequalities with nonnegative weights $(1-\beta)$ and $\beta$, respectively, and then subtracting the second inequality, we get $2 \bar{\gamma} \leq w^\top (\beta x + (1-\beta)v_+ - v_-) \leq \|w\|_* \left\| \beta x + (1-\beta)v_+ - v_- \right\|$.
Since $\|w\|_* = 1$ and the preceding relation holds for all $\beta \in [0,1]$, we deduce that $2 \bar{\gamma}$ is less than or equal to
$\min_{\beta \in [0,1]} \left\| \beta x + (1-\beta)v_+ - v_- \right\| = \tau(\{v_+,x\}, \{v_-\}) = 2 \gamma(\{v_+,x\}, \{v_-\})$. Now recall that $\|v_+ - v_-\| = 2\gamma>0$, so we can rewrite $\min_{\beta \in [0,1]} \left\| \beta x + (1-\beta)v_+ - v_- \right\| = 2 \gamma \min_{\beta \in [0,1]}  \left\| \frac{v_+ - v_-}{\|v_+ - v_-\|} - \frac{\beta (v_+ - x)}{\|v_+ - v_-\|} \right\| \leq 2\gamma  \kappa_\circ\left( \frac{\|v_+ - x\|}{\|v_+ - v_-\|}, \frac{\ell(v_+ - v_-)^\top(v_+ - x)}{\|v_+ - v_-\|} \right)$,
where the last inequality follows from the definition of $\kappa_{\circ}$. An analogous proof follows when $x$ is added to $V_-$ instead of $V_+$.
\Halmos
\endproof

The crux of our analysis is to show that $\kappa_\circ$ is, in fact, $< 1$ when its inputs are $> 0$.
\begin{lemma}\label{lemma:kappa}
For $\delta,\eta \in \bbR$,
\begin{align*}
\kappa_{\circ}(\delta,\eta) \begin{cases}
= -\infty, & \text{if } \delta < \eta,\\
< 1, &  \text{if } \delta \geq \eta > 0,\\
= 1, &  \text{if } \eta \leq 0.
\end{cases}
\end{align*}
\end{lemma}
\proof{Proof of \cref{lemma:kappa}.}
When $\ell(u)^\top z \geq \eta$ and $\|z\| \leq \delta$, by H\"{o}lder's inequality we have $\eta \leq \|\ell(u)\|_* \|z\| = \|z\| \leq \delta$ (recall that by definition of $\ell(u)=\nabla\|u\|$ and convex analysis we always have $\|\ell(u)\|_*=1$). So, the maximum defined in the definition of $\kappa_{\circ}(\delta,\eta)$ is infeasible when $\delta < \eta$. When $\eta \leq 0$, we can choose $z = 0$ to make the minimum in $\kappa_{\circ}(\delta,\eta)$ definition 
equal to $1$, resulting in $\kappa_{\circ}(\delta,\eta) = 1$.

Now consider the case $\delta \geq \eta > 0$. Fix $u,z$ such that $\|u\| = 1$, $\|z\| \leq \delta$, and $\ell(u)^\top z \geq \eta$. Define $f(\beta) := \left\| u - \beta z \right\|$.
Clearly, $f(\beta)$ is a differentiable convex univariate function and $f(0) = \|u\| = 1$. We will show that $f'(0) < 0$, which then shows that the minimum in $\kappa_{\circ}(\delta,\eta)$ definition must be $< 1$ since we can increase $\beta$ by a small amount to get $f(\beta) < 1$. Note that $f'(0)$ can be computed via the directional derivative formula: $f'(0) = -\ell(u)^\top z \leq -\eta < 0$. Therefore, $\min_{\beta \in [0,1]} \|u-\beta z\| < 1$.
Now since $\|\cdot\|$ is continuous, $\min_{\beta \in [0,1]} \|u-\beta z\|$ is a continuous function of $\beta$. Furthermore, the domain $\{(u,z) :~ \|u\|=1,\, \|z\| \leq \delta,\, \ell(u)^\top z \geq \eta \}$ is compact since $\ell(u)^\top z$ is continuous in $(u,z)$. Therefore, maximizing $\min_{\beta \in [0,1]} \|u-\beta z\|$ over this domain is solvable at some $(u^*,z^*)$. Moreover, $\kappa_{\circ}(\delta,\eta) = \min_{\beta \in [0,1]} \|u^*-\beta z^*\| < 1$ as well.
\Halmos
\endproof

We now provide bounds on $m(\cdot)$ (which also yields an upper bound on the number of mistakes). To ease notation, recall that in \cref{ass:margin-diameter} we defined $r := D/\gamma_*$, and we now also define
\begin{align}\label{eq:kappa}
\kappa(r,\phi) := \kappa_\circ\left( r/2, (1-\phi)/2 \right).
\end{align}
As a consequence of the previous lemmas, we have the following relationship between $\gamma_t,\gamma_{t+1}$ throughout our algorithm, and importantly, this holds \emph{regardless} of whether we use the na\"{i}ve or efficient implementation in \cref{alg:max-margin}, which is a key step to ensuring computational efficiency.
\begin{lemma}\label{lemma:max-margin-increment-alg}
Suppose that \cref{alg:max-margin} is run with aggressiveness parameter $\rho \in [0,1]$, and recall that $a_t := y_t (w_t^\top x_t + b_t)$. Let $\phi \leq \rho$ be a fixed constant.
For $t \geq 3$, whenever $a_t < \phi \gamma_t \leq \rho \gamma_t$, we have
\[ \frac{\gamma_{t+1}}{\gamma_t} \leq \kappa_{\circ}\left( \frac{D}{2\gamma_t}, \frac{\gamma_t - a_t}{2\gamma_t} \right) \leq \kappa(r,\phi). \]
\end{lemma}
\proof{Proof of \cref{lemma:max-margin-increment-alg}.}
We prove the first inequality. For $t \geq 3$, at time $t-1$, we have $v_{t,+} \in \conv(V_{t,+}),v_{t,-} \in \conv(V_{t,-})$. 
Since $a_t < \phi \gamma_t \leq \rho \gamma_t$, by the definition of \cref{alg:max-margin}, $(w_{t+1},b_{t+1})$ is obtained by solving \cref{eq:margin} with $V_{t+1,+},V_{t+1,-}$. Let $k \leq t-1$ be the last prior iteration with an update, i.e., $a_k < \rho \gamma_k$. So all data at iteration $k+1$ are identical to time $t$. In particular, no points are added to $V_{s,+},V_{s,-}$ sets from $s=k+1,\ldots,t$.

For the na\"{i}ve implementation of \cref{alg:max-margin}, \cref{lemma:max-margin-increment} can be applied immediately, since there is only one point $x_t$ difference between the sets $V_{t+1,+},V_{t+1,-}$ used to compute $(w_{t+1},b_{t+1})$ and the previous sets $V_{k+1,+},V_{k+1,-}$ used to compute $(w_{k+1},b_{k+1}) = (w_t,b_t)$. 
Furthermore, $\gamma_t = \gamma_{k+1} = \gamma(V_{k+1,+}, V_{k+1,-})$. For the efficient implementation, notice that there is one point difference between $V_{t+1,+} \cup V_{t+1,-}$ and $\{v_{k+1,+},v_{k+1,-}\} = \{v_{t,+},v_{t,-}\}$, and $\gamma_t = \gamma_{k+1} = \gamma(\{v_{k+1,+}\}, \{v_{k+1,-}\}) = \gamma(\{v_{t,+}\}, \{v_{t,-}\})$. So \cref{lemma:max-margin-increment} can also be applied. Noting that in both cases $w_t = \ell(v_{t,+} - v_{t,-})$ and $2\gamma_t = \|v_{t,+} - v_{t,-}\|$, this gives $\frac{\gamma(V_{t+1,+},V_{t+1,-})}{\gamma_t} \leq
\kappa_\circ\left( \frac{\|v_{t,+} - x_t\|}{2\gamma_t}, \frac{w_t^\top (v_{t,+} - x_t)}{2\gamma_t} \right)$ when $y_t = +1$, and $\frac{\gamma(V_{t+1,+},V_{t+1,-})}{\gamma_t} \leq \kappa_\circ\left( \frac{\|x_t - v_{t,-}\|}{2\gamma_t}, \frac{w_t^\top (x_t-v_{t,-})}{2\gamma_t} \right)$ when $y_t = -1$.
Notice now that in \cref{alg:max-margin} $v_{t,+},v_{t,-}$ are always convex combinations of positive and negative points seen so far, thus $\|v_{t+1,+} - x_t\| \leq D$ when $y_t = +1$ and $\|x_t - v_{t+1,-}\| \leq D$ when $y_t = -1$. Furthermore, when $y_t = +1$, we have $w_t^\top v_{t,+} + b_t = \gamma_t$ and $w_t^\top x_t + b_t = a_t$, therefore $w_t^\top (v_{t,+} - x_t) = \gamma_t - a_t$.
Similarly, when $y_t = -1$, we have $w_t^\top v_{t,-} + b_t = -\gamma_t$ and
$w_t^\top x_t + b_t = -a_t$, therefore $w_t^\top (x_t - v_{t,-}) = \gamma_t - a_t$.
Since $\kappa_\circ$ is non-decreasing in the first argument, we have for both na\"{i}ve and efficient implementations,
$\gamma_{t+1}/\gamma_t = \frac{\gamma(V_{t+1,+}, V_{t+1,-})}{\gamma_t} \leq \kappa_{\circ} \left( \frac{D}{2\gamma_t}, \frac{\gamma_t-a_t}{2\gamma_t}  \right)$.

To obtain the second inequality, observe that since $a_t \leq \phi \gamma_t$, we have $(\gamma_t - a_t)/(2\gamma_t) \geq (1-\phi)/2 > 0$. Furthermore, by \cref{lemma:gamma_bound} we have $\gamma_t \geq \gamma_*$ and $\gamma_*>0$ as the data are separable, therefore $r = D/\gamma_* \geq D/\gamma_t$. Then, since $\kappa_{\circ}$ is non-decreasing in the first argument and non-increasing in the second argument, we deduce that $\kappa_{\circ}\left( \frac{D}{2\gamma_t}, \frac{\gamma_t - a_t}{2\gamma_t} \right) \leq \kappa_{\circ}(r/2, (1-\phi)/2) = \kappa(r,\phi)$.
\Halmos
\endproof

As a consequence of the previous lemmas, we arrive at the following mistake bound for \cref{alg:max-margin}.
\begin{theorem}\label{thm:margin-performance}
Fix $\rho \in [0,1]$. Then, \cref{alg:max-margin} with aggressiveness parameter $\rho$ \emph{simultaneously} satisfies for all $\phi \leq \rho$
\[ m(\phi) \leq 
\begin{cases} 0, & \text{if } \phi < 1-r , \\
\frac{\log(\gamma_3/\gamma_*)}{-\log(\kappa(r,\phi))}, & \text{if } \phi \geq 1-r.
\end{cases} \]
\end{theorem}
\proof{Proof of \cref{thm:margin-performance}.}
Recall that $a_t := y_t (w_t^\top x_t + b_t)$.
Note that by \cref{lemma:kappa}, $\kappa_{\circ}(r/2, (1-\phi)/2) = -\infty$ when $\phi < 1-r$, so in this case $m(\phi) = 0$. 
When $\phi \geq 1-r$ and $a_t \leq \phi \gamma_t$ we have $\gamma_{t+1}/\gamma_t \leq \kappa(r,\phi)$ from \cref{lemma:max-margin-increment-alg}, and $\gamma_{t+1}/\gamma_t \leq 1$ otherwise. Taking the logarithm of both sides and summing from $t=3,\ldots,\infty$, we get
$ \lim_{t \to \infty} \log(\gamma_{t+1})  - \log(\gamma_3) \leq m(\phi) \log(\kappa(r,\phi)). 
$ As $\lim_{t \to \infty} \log(\gamma_{t+1}) \geq \log(\gamma_*)$, the result then follows after rearrangement.
\Halmos
\endproof

\section{Computational considerations and theoretical complexity comparisons}\label{sec:implementation}

In this section, we examine the case when the norm $\|\cdot\|$ is chosen to be the $\ell_2$-norm, describe how to implement the efficient version of \cref{alg:max-margin},
and provide an explicit performance guarantee based on parameters in \cref{ass:margin-diameter} based on the $\ell_2$-norm. We then provide a comparison of our guarantees against the existing algorithms and their associated guarantees.

\subsection{Implementation details and guarantees for \texorpdfstring{$\ell_2$}{l2}-norm} \label{sec:implementation-l2}

Each iteration of \cref{alg:max-margin} requires us to solve \cref{eq:margin}, or equivalently solve \cref{eq:convex-hull-distance}.
When the norm is $\ell_p$ for $p\in(1,\infty)$, \cref{eq:convex-hull-distance} can be cast as a power cone optimization problem, thus both na\"{i}ve and efficient implementations can be solved using off-the-shelf conic solvers. However, the size of $V_+ \cup V_-$ for the efficient implementation is only three points, thus \cref{eq:convex-hull-distance} is equivalent to solving the problem
\begin{align}\label{eq:convex-hull-simple}
\min_{\beta \in [0,1]} \|u - \beta z\|_p
\end{align}
for given vectors $u,z$. Indeed, we will see that we can even avoid using conic solvers by analyzing \cref{eq:convex-hull-simple} further. In particular, we are able to provide a closed-form solution for \eqref{eq:convex-hull-simple} when the norm is the $\ell_2$-norm. If a general $\ell_p$-norm is used for $p \in (1,\infty)$, then \cref{eq:convex-hull-distance} can also be efficiently solved via one-dimensional line search (for brevity, we omit these details). Furthermore, since \cref{eq:convex-hull-simple} also appears in the definition of $\kappa(\cdot)$ in \cref{eq:kappa}, our analysis allows us to derive an explicit form for the guarantee in \cref{thm:margin-performance} in terms of the parameters $r,\gamma_*$.

When $p=2$, the minimizer $\beta_*(u,z)$ of \cref{eq:convex-hull-simple} must satisfy the optimality condition $(\beta_*(u,z)\|z\|_2^2 - u^\top z)(\beta - \beta_*(u,z)) \geq 0$ for all $\beta \in [0,1]$, which results in the closed-form solution
\[\beta_*(u,z) = 
\begin{cases}
    0, & \text{if } u^\top z < 0, \\
    \frac{u^\top z}{\|z\|_2^2}, & \text{if } 0 \leq u^\top z \leq \|z\|_2^2, \\
    1, & \text{if } u^\top z > \|z\|_2^2.
\end{cases} \]
Further to \cref{rem:complexity}, we see that solving \cref{eq:convex-hull-distance} in the efficient implementation of \cref{alg:max-margin} has an operation complexity of $O(d)$.

The closed form for $\beta_*(u,z)$ also implies that
\[ \cref{eq:convex-hull-simple} = \begin{cases} 
\|u\|_2, & \text{if } u^\top z < 0,\\
\sqrt{\|u\|_2^2 - \frac{(u^\top z)^2}{\|z\|_2^2}}, & \text{if } 0 \leq u^\top z \leq \|z\|_2^2,\\ 
\|u-z\|_2, & \text{if } u^\top z > \|z\|_2^2. 
\end{cases} \]
Recall that $\kappa_\circ(\delta,\eta)$ chooses $u,z$ to maximize the optimal value of \cref{eq:convex-hull-simple}, subject to $\|u\|_2 = 1$, $\|z\|_2 \leq \delta$, $\ell(u)^\top z = u^\top z \geq \eta$. When $\delta< 0$ or $\delta < \eta$, due to the Cauchy-Schwarz inequality the problem for computing $\kappa_\circ(\delta,\eta)$ is infeasible, hence $\kappa_\circ(\delta,\eta) = -\infty$ in these cases. When $\eta < 0$, note that the minimum is always $\|u\|_2 = 1$, so $\kappa_\circ(\delta,\eta) = 1$ in this case. In other cases (i.e., when $0 \leq \eta \leq \delta$), note that we can always choose $u,z$ such that $\|u\|_2=1$, $\|z\|_2 = \delta$, $u^\top z = \eta$, simply by setting $u = e_1$, and $z = \eta e_1 + \sqrt{\delta^2 - \eta^2} e_2$, where $e_1,e_2$ are the first two standard basis vectors. If, additionally, we have $\eta \leq \delta^2$, then any $u,z$ satisfying the constraints will yield $\cref{eq:convex-hull-simple} = \sqrt{\|u\|_2^2 - (u^\top z)^2/\|z\|_2^2} \leq \sqrt{1 - \eta^2/\delta^2}$. Therefore $\kappa_\circ(\eta,\delta) = \sqrt{1 - \eta^2/\delta^2}$ as there exist $u,z$ to make the constraints tight. Finally, if additionally we have $\delta^2 < \eta \leq \delta$, then $\cref{eq:convex-hull-simple} = \|u-z\|_2 = \sqrt{\|u\|_2^2 - 2u^\top z + \|z\|_2^2} \leq \sqrt{1 - 2\eta + \delta^2}$. Again, since there exist $u,z$ to make the constraints tight, we deduce that $\kappa_\circ(\delta,\eta) = \sqrt{1 - 2\eta + \delta^2}$ in this case. In summary, we have
\begin{align*}
    \kappa_\circ(\delta, \eta)
    = \begin{cases}
    1, & \text{if } \eta < 0 \leq \delta, \\
    \sqrt{1-\frac{\eta^2}{\delta^2}}, & \text{if } 0 \leq \eta \leq \min\{\delta^2,\delta\}, \\
    \sqrt{1+\delta^2-2\eta},  & \text{if } \delta^2< \eta \leq \delta,\\
    -\infty, & \text{if } \max\{0,\eta\} > \delta.
\end{cases}
\end{align*}
Substituting $\delta = \frac{r}{2}$ and $\eta = \frac{1-\phi}{2}$ into $\kappa(r,\phi) = \kappa_\circ \left(\frac{r}{2}, \frac{1-\phi}{2} \right)$, where $r \geq 0$ and $\phi < 1$ by assumption, and then substituting the resulting expression of $\kappa_\circ \left(\frac{r}{2}, \frac{1-\phi}{2} \right)$ from the preceding formula into \cref{thm:margin-performance}, we obtain the following bound on $m(\phi)$ in the case of $\ell_2$-norm:
\begin{align*}
    m(\phi) &\leq  
    \begin{cases}
    0, &  \text{if } \frac{r}{1-\phi} < 1,\\
    \frac{2\log(\gamma_3/\gamma_*)}{-\log\left(\phi+\frac{r^2}{4} \right)}, & \text{if } 1 \leq \frac{r}{1-\phi} < \sqrt{2},\\
    \frac{2\log(\gamma_3/\gamma_*)}{-\log\left(1- \left(\frac{1-\phi}{r}\right)^2 \right)}, & \text{if } \frac{r}{1-\phi} \geq \sqrt{2}.
    \end{cases}
\end{align*}

\subsection{Comparison of theoretical guarantees}\label{sec:comparison}

In \cref{tab:perf_guarantees}, for the case of $\ell_2$-norm, we provide a comparison of theoretical guarantees (in terms of the \# of mistakes and the maximum \# of iterations needed to ensure a margin of $\phi\gamma_*$) for the efficient (\eomm) and na\"{i}ve (\nomm) versions of \cref{alg:max-margin} along with five existing online classification algorithms: the classical \pctr algorithm \citep{rosenblatt_perceptron_1958}, \romma, aggressive ROMMA (\aromma) \citep{li_relaxed_1999}, \pumma \citep{IshibashiEtAl2008} and \alma \citep{gentile_new_2000}. We also report in this table whether the corresponding algorithm requires and uses the a priori knowledge of $b_*=0$ or not.

\newcommand*{\mline}[1]{\begingroup
    \renewcommand*{\arraystretch}{1.1}   \begin{tabular}[c]{@{}>{\raggedright\arraybackslash}p{2cm}@{}}#1\end{tabular}  \endgroup
}
\begin{table}[htp]
    \centering
    \resizebox{\textwidth}{!}{
    \begin{tabular}{|l|l|l|l|}
    \hline
        & Mistake bound & Max \# iters to reach at least $\phi \gamma_*$ margin         & $b_*$\\
    \hline
    \begin{tabular}{c}
    (Alg. \ref{alg:max-margin})\\
    \nomm\\
    \& \eomm
    \end{tabular} & $
    \begin{cases}
    2, & \text{if } r < 1,\\
        2+\frac{\log(\gamma_3/\gamma_*)}{-\log\left(r/2 \right)}, & \text{if } 1 \leq r < \sqrt{2},\\
        2+\frac{2\log(\gamma_3/\gamma_*)}{-\log\left(1- \frac{1}{r^2} \right)}, & \text{if } r \geq \sqrt{2}.
    \end{cases}$ & $
    \begin{cases}
        2, & \text{if } \frac{r}{1-\phi} < 1-\phi,\\
        2+\frac{2\log(\gamma_3/\gamma_*)}{-\log\left(\phi+\frac{r^2}{4} \right)}, & \text{if } 1 \leq \frac{r}{1-\phi} < \sqrt{2},\\
        2+\frac{2\log(\gamma_3/\gamma_*)}{-\log\left(1- \left(\frac{1-\phi}{r}\right)^2 \right)}, & \text{if } \frac{r}{1-\phi} \geq \sqrt{2}.
    \end{cases}$ & free \\
    \hline
    \pumma & \multirow{4}{*}{$\displaystyle O\left( \frac{\barD^2}{\gamma_*^2} \right) = O \left(r^2 \left(\frac{\barD}{D}\right)^2 \right) $} & \multirow{2}{*}{$\displaystyle O \left(\left( \frac{r}{1-\phi} \cdot \frac{\barD}{D} \right)^2 \right)$} & free\\
    \cline{4-4}
    \aromma & & & \\
    \cline{1-1} \cline{3-3}
    \romma &  & \multirow{2}{*}{None} & \multirow{2}{*}{$0$} \\
    \pctr &  &  & \\
    \cline{1-3}
    \alma & $\displaystyle O\left(\frac{1}{(\hat{\gamma}_*)^2}\right)$ for normalized margin $\hat{\gamma}_*$ & $\displaystyle O \left(\frac{1}{(\hat{\gamma}_*)^2(1-\phi)^2}\right)$ & \\
    \hline
    \end{tabular}
    }
    \caption{Theoretical guarantees for various algorithms in the case of $\ell_2$-norm. \alma refers to \alma($1-\phi;\sqrt{8}/(1-\phi),\sqrt{2}$) a la \citet[Theorem 3]{gentile_new_2000}, $\gamma_*$ is the maximum margin on the original data while $\hat{\gamma}_*$ is the max-margin on the normalized data $x_t/\|x_t\|_2$.         } 
    \label{tab:perf_guarantees}
\end{table}

As mentioned in \cref{rem:translation-invariance,rem:translation-invariance-alg}, \cref{alg:max-margin} is translation invariant, and consequently the bounds on $m(\phi)$ depend only on $\gamma_3$, $\gamma_*$ and $r$, which are all translation invariant quantities. In contrast, the guarantees given for other algorithms (except \alma) depend on the ratio $\barD/D$, which is a main point of difference between our results and the previous bounds. Since $\barD$ is not a translation invariant quantity, $\barD/D$ can be thought of as measuring how much the lack of translation invariance may affect these other algorithms' performances. 
The exception is \alma, which works with normalized feature vectors $x_t/\|x_t\|$ and aims to find the maximum margin $\hat{\gamma}_*$ on the normalized data. (Note that $\hat{\gamma}_* \neq \gamma_*$ in general.) This data remains separable if the original maximum margin classifier has $b_* = 0$. However, it is easy to see that \alma is also not translation invariant, since normalizing a vector is not translation invariant. We provide a detailed comparison of the mistake bound of \cref{alg:max-margin} against others in \cref{sec:EC-guarantees}.

This translation invariance (or lack thereof) difference between our mistake bounds versus the ones from the literature has significant implications. The theoretical gap in mistake bounds between the $r^2 (\barD/D)^2$ bound of \pctr/\romma/\pumma-type algorithms and those of \cref{alg:max-margin} can be made arbitrarily large by constructing examples in which $D$ and $\gamma_*$ remain fixed while $\barD$ grows without bound. We now present such an example, which also illustrates that the arbitrarily large discrepancy in the mistake bounds is not merely a theoretical possibility, but can also arise in practice. Due to space constraints, we consider only the \pctr and omit the consideration of other algorithms. 

\begin{example}\label{ex:bad-perceptron}
Fix parameters $c,r > 0$ to be chosen later. Suppose that there are only three possible feature-label pairs: $z^1 = (x^1,y^1) = ((c,1),+1)$, $z^2 = (x^2,y^2) = ((c,-1),-1)$ and $z^3 = (x^3,y^3)=((c+r,-1),-1)$. Then $(w_*,b_*) = ((0,1),0)$, $\gamma_* = 1$, $\barD = \sqrt{(c+r)^2+1}$, and $D = D/\gamma_* = r$.

We will set $r = 2/c$, and $c > 2$. Since $r < 1$, both versions of \cref{alg:max-margin}, with any aggressiveness parameter $\rho \in [0,1]$, will make at most two mistakes on \emph{any} sequence where each $(x_t,y_t) \in \{z^1,z^2,z^3\}$. (In fact, for any $c,r > 0$ \cref{alg:max-margin} commits at most three mistakes, but for brevity, we omit this argument.)

In contrast, we now show that for any $m \in \bbN$, we can construct a sequence for which the \pctr algorithm makes $\geq m$ mistakes, thus the gap in terms of the mistake bounds is arbitrarily large. The \pctr algorithm (using the a priori information of $b_*=0$) proceeds as follows: initialize $w_1 = (0,0)$, then update $w_{t+1} = w_t + y_t x_t$ only when $y_t w_t^\top x_t \leq 0$. We assume that the sequence $(x_t,y_t)$ alternates between $z^1$ and $z^3$ until some time $m$ to be defined later. More precisely, $x_t = (c,1)$, $y_t=+1$ for $t \leq m$ odd, and $x_t = (c+r,-1)$, $y_t = -1$ for $t \leq m$ even. We will choose the parameters $c,r$ to ensure that mistakes are made for all $t \leq m$, thus the number of mistakes is $\geq m$. To do this, it is easy to show that if mistakes occur, we have $w_t = (c-(t/2-1)r,t-1)$ when $t$ is even, and $w_t = (t-1)(-r/2,1)$ when $t$ is odd. Therefore mistakes will occur when $(c-(t/2-1)r,t-1)^\top (c+r,-1) \geq 0 \iff t \leq \frac{2(c^2 + r^2 + 2cr + 1)}{cr + r^2 + 2}$ and $(t-1)(-r/2,1)^\top (c,1) \leq 0 \iff cr \geq 2$. The second condition is guaranteed since $r = 2/c$. We then set $m := \left\lfloor \frac{c^2 + 4/c^2 + 5}{2(1 + 1/c^2)} \right\rfloor$, and thus the first condition is also satisfied and mistakes occur whenever $t \leq m$. Therefore, we can force the \pctr algorithm to make $\geq m = \Omega(c^2)$ mistakes. Note that the mistake bound is $r^2 (\barD/D)^2 = (c+r)^2 + 1 = c^2 + 4/c^2 + 3 = O(c^2)$, thus this example shows that the mistake bound for the \pctr is asymptotically tight.
\epr
\end{example}

In \cref{sec:numerical} we numerically investigate the effect that $\barD/D$ and non-zero $b_*$ have on the performance of different algorithms on real data. We note that the results in \cref{sec:numerical} suggest that \pumma is also translation invariant; inspired by these observations, we verify that this is indeed the case in \cref{sec:translation-invariance-pumma}. We also conjecture that the $r^2 (\barD/D)^2$ mistake bound of \pumma can be improved to be based on only translation invariant quantities.  

\section{Numerical study}\label{sec:numerical}

We conducted a numerical study comparing our proposed algorithms, \eomm and \nomm as well as \ceomm which is a conservative version of \eomm with $\rho=0$, against five benchmark methods listed in \cref{tab:perf_guarantees} using real-world data and taking $\|\cdot\|$ to be the $\ell_2$-norm. We provide a high-level description of our study here; see the e-companion \cref{sec:exp-setup} for full implementation details (e.g., on hyper-parameter tuning and adaptations to handle the $b_* \neq 0$ case for other algorithms). 

We use the Adult dataset from the UCI Machine Learning Repository \citep{adult_dataset}, with each data point $(x_t,y_t)$ representing a feature-label pair in an online classification task. After preprocessing to remove missing values, encode categorical features, and enforce linear separability (see  \cref{sec:exp-preprocessing} for full details of this preprocessing), we obtain a dataset of $35508$ points in $d = 96$ dimensional space, where 27.11\% of the points have label +1.

To analyze algorithm performance under varied conditions, we construct dataset variants by applying the following three transformations (see \cref{sec:exp-transformations} for full descriptions of these transformations): 
\textbf{(i) Margin normalization:} We transform the data to ensure that the true margin is $\gamma_*=1$ without changing the diameter bound $D$.
\textbf{(ii) Bias control:} We either retain the original bias ($b_*\neq 0$) of the best margin classifier $(w_*,b_*)$  or transform the data to enforce that $b_*=0$.
\textbf{(iii) Feature scaling:} For an input parameter of $\theta \in \{0,0.25,0.5,0.75,1\}$, we amplify $\barD$, while leaving $D$, $r$, and $\gamma_*$ unchanged. This tests sensitivity to the $\barD/D$ ratio, a key term in terms of theoretical performance guarantees of \cref{alg:max-margin} and the existing methods (see \cref{tab:perf_guarantees}). 
Combining five values of $\theta$ with both bias options gives us ten different variants of our dataset. Margin normalization to ensure $\gamma_* = 1$ is performed in all variants of the dataset. 
Note that $D$ is invariant in all of these transformations, and as reference, $D=220.43$ in our dataset.

As computational performance measures, we record the following metrics at the end of each algorithm's run: (i) the margin of the final classifier on the entire dataset, i.e., $\bar{\gamma} := \min_{x,y} \gamma(x,y; w_T,b_T)$, where $(w_T,b_T)$ is the last classifier returned by the algorithm and the minimum is taken over all points $(x,y)$ in the dataset; (ii) the total number of mistakes (misclassifications) $m$ the algorithm makes throughout its run; (iii) the number of iterations $\tau$ it takes to reach a classifier with positive margin $\bar{\gamma} > 0$; and (iv) the total computation time in seconds.
Results for $\theta = 0$ over a single pass of the dataset are reported in \cref{tab:perf-single-pass}, while plots of the first three metrics for various $\theta$ are shown in \cref{fig:plot-metrics}. We do not plot the run times, since they do not vary much for different $\theta$ or $b_*$ values. The full table of metrics (including the run times) for all $\theta$ values is provided in \cref{sec:extra_numerical_results}.

\begin{table}[htbp]
    \centering
    \begin{tabular*}{\linewidth}{@{\extracolsep{\fill}}l|rrrr||rrrr}
    \toprule
        & \multicolumn{4}{c}{$b_* = 0$} & \multicolumn{4}{c}{$b_* = -0.2848$}\\
    \cline{2-5} \cline{6-9} 
        & $m$ & $\bar{\gamma}$ & $\tau$ & time (s) & $m$ & $\bar{\gamma}$ & $\tau$ & time (s) \\
    \midrule
    $\theta = 0$ & \multicolumn{4}{c}{$\barD/D = 0.9666$} & \multicolumn{4}{c}{$\barD/D = 0.9658$} \\
    \hline
    \nomm & 4 & 1.00 & 70 & 76.56 & 4 & 1.00 & 70 & 74.25 \\
    \ceomm & 21 & 0.05 & 16808 & 0.10 & 21 & 0.05 & 16808 & 0.11 \\
    \eomm & 6 & 0.84 & 534 & 0.11 & 6 & 0.84 & 534 & 0.12 \\
    \romma & 38 & 0.16 & 28065 & 0.09 & 4114 & -- & -- & 0.22 \\
    \aromma & 11 & 0.73 & 265 & 0.15 & 2643 & -- & -- & 0.21 \\
    \pumma & 9 & 0.68 & 584 & 0.08 & 9 & 0.68 & 584 & 0.08 \\
    \pctr & 59 & -- & -- & 0.07 & 52 & -- & -- & 0.07 \\
    \alma & 53 & -- & -- & 0.07 & 174 & -- & -- & 0.07 \\
    \bottomrule
    \end{tabular*}
    \caption{
    Numerical results for $\theta=0$ when the algorithms are allowed a single pass over the dataset. A dash for $\bar{\gamma}$ indicates that the final classifier does \emph{not} have a positive margin on the dataset; similarly, a dash for $\tau$ indicates that no classifier generated by that algorithm throughout the run has positive margin.
    }
    \label{tab:perf-single-pass}
\end{table}

\begin{figure}[tp]
        \centering
        \includegraphics[width=.95\textwidth]{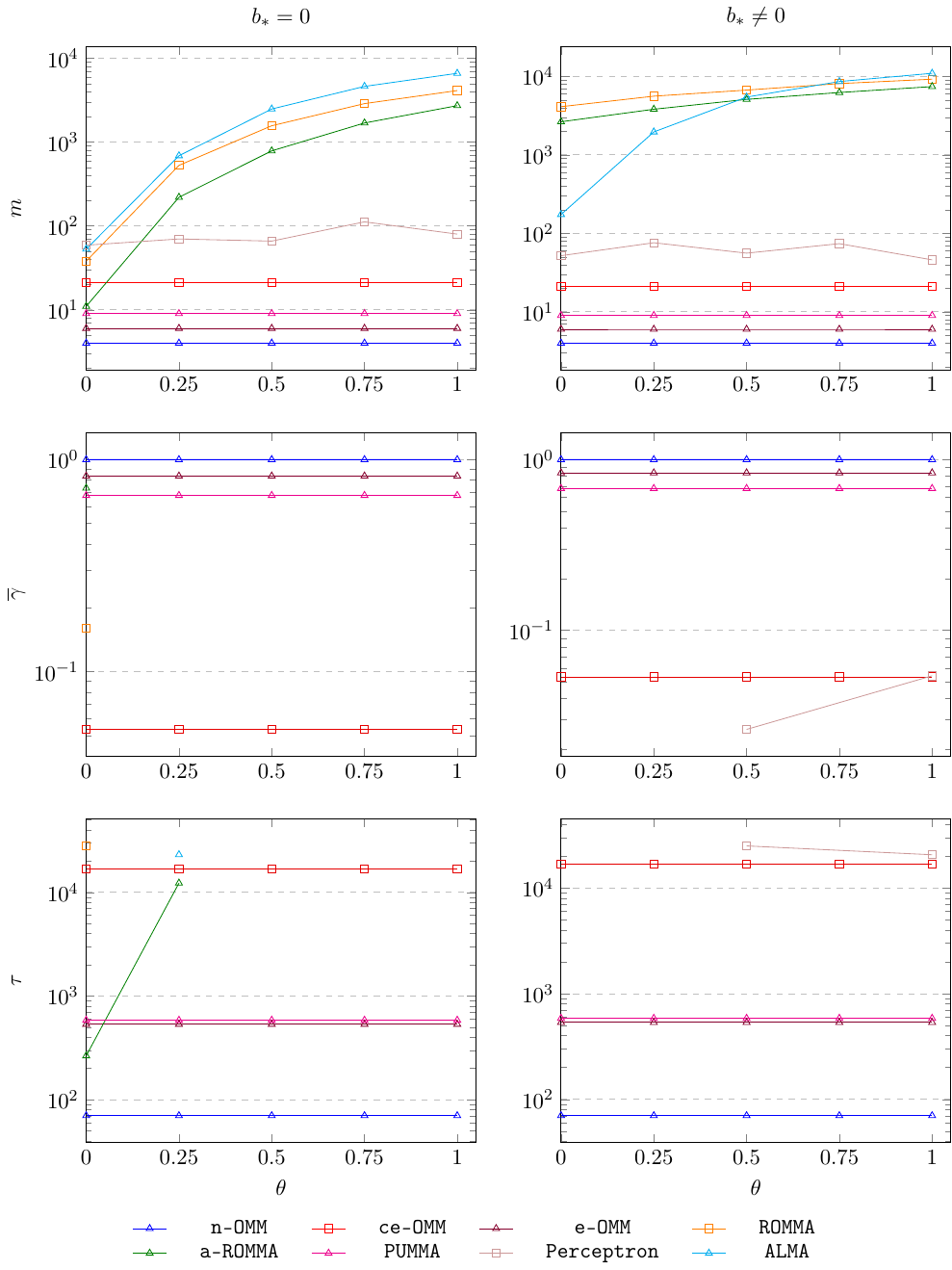}
        \caption{Plots of the number of mistakes $m$, the margin $\bar{\gamma}$, and the number of iterations to needed find a classifier with positive margin $\tau$ (if any) within a single run over the data set variants generated by different values of $\theta$ and whether $b_*=0$ or not.
    }\label{fig:plot-metrics}
\end{figure}

We observe that the runtime of \nomm at around $\sim$70s is significantly higher than other algorithms, which finish a pass in less than 0.5 seconds. On the other hand, \nomm makes only four mistakes and recovers the perfect margin $\gamma_* = 1$. This is expected since on its last iteration, it is essentially solving an offline maximum margin problem on the whole dataset. Examining \cref{fig:plot-metrics}, we notice that the metrics of \nomm, \ceomm, \eomm, and \pumma remain constant as we change $\theta$ and $b_*$. This is expected for \nomm, \ceomm, and \eomm by \cref{rem:translation-invariance-alg}, since these transformations are only translations of the data $x_t$ by a fixed vector. This behavior suggests that \pumma is also translation invariant; inspired by this observation, we verify that this is indeed the case in \cref{sec:translation-invariance-pumma}.
Among the translation-invariant algorithms other than \nomm, we note that \eomm performs best with 6 mistakes and the recovered margin of 0.84, followed by \pumma with 9 mistakes and a margin of 0.68, and then by the conservative version of our algorithm \ceomm with 21 mistakes.

Comparing now the translation-invariant algorithms (\nomm, \ceomm, \eomm, \pumma) against the other algorithms (\romma, \aromma, \pctr, \alma), we notice that the translation-invariant algorithms consistently outperform the others in all metrics ($m$, $\bar{\gamma}$, and $\tau$) even under the most favorable setting of $\theta=0$. There are two notable exceptions: \aromma achieves a higher margin than \pumma, while both \romma and \aromma achieve higher margin than \ceomm in the most favorable case for non-translation invariant algorithms when $\theta=0$, $b_* = 0$. Note that \ceomm and \romma are \emph{conservative} algorithms, updating only when there is a mistake. Therefore, the margin is not a meaningful measure here, as the algorithms will stop updating once all data points are correctly classified, i.e., $\bar{\gamma} > 0$. For conservative algorithms such as \ceomm, \romma and \pctr, the number of iterations $\tau$ until $\bar{\gamma} > 0$ is a more meaningful measure, and \ceomm is lower than other conservative algorithms. When $\theta$ increases, or $b_* \neq 0$, we notice in \cref{fig:plot-metrics} that \romma, \aromma and \alma degrade sharply, performing even worse than the conservative \ceomm.
Overall, when $\theta>0$ or $b_*\neq0$, the number of mistakes made by \romma, \aromma and \alma exceed 1000 or 10,000 quickly and they are unable to recover a correct classifier in these setting in a single pass over the dataset.
The number of mistakes for the \pctr remains surprisingly stable, though note that it is worse than \ceomm, and it fails to find a correct classifier within one pass over the data, i.e., $\tau$ is undefined.

In an additional experiment, we test the impact of going through the same dataset five times, each with a new random permutation of the points, except \nomm where only a single pass through the dataset is made, as it will stop updating after seeing all the data once. We report our findings in \cref{sec:extra_numerical_results}. 
The findings are broadly consistent with the one-pass experiments: the translation invariant algorithms \eomm and \pumma generally perform well, but the other algorithms either do not perform well to begin with (\alma, \pctr) and/or severely degrade as $\theta$ increases and $b_* \neq 0$ (\romma, \aromma). See \cref{sec:extra_numerical_results} for a full discussion. 

Collectively, these observations highlight three findings. First, the importance of translation invariance in number of mistakes and margin recovery is demonstrated through the superior performance of \nomm, \ceomm, \eomm and \pumma. Second, although conservative algorithms such as \ceomm directly target the number of mistakes, using aggressive ones such as \eomm can result in lower number of mistakes. Third, while both algorithms are translation invariant, our \eomm algorithm outperforms \pumma (more significantly in terms of margin recovery) by making judicious use of the dual certificates to accurately summarize past information on the margin.

\newpage
\appendix

\section{Further comments on \texorpdfstring{\cref{tab:perf_guarantees}}{Table 1}}\label{sec:EC-guarantees}

Let us qualitatively compare the $O(r^2 (\barD/D)^2)$ mistake bound of \pctr/\romma/\pumma algorithms with the mistake bound of \cref{alg:max-margin}, summarized in \cref{tab:perf_guarantees}.

When $r < 1$, i.e., when ${D\over \gamma_*} <1$, \cref{alg:max-margin} makes at most two mistakes, whereas the bound $O(r^2 (\barD/D)^2)$ can be arbitrarily high depending on how large $\barD/D$ is. \cref{ex:bad-perceptron} has shown that we can construct examples for $r < 1$ where $\barD/D$ is arbitrarily large, and the perceptron demonstrably makes $\Omega(r^2(\barD/D)^2)$ mistakes, i.e., this bound is tight for the perceptron in the $r < 1$ regime.

When $r \geq 1$, our bound has an additional $\log(\gamma_3/\gamma_*)$ term. It can be seen that $\gamma_3 \leq 2(D+\gamma_*)$, hence $\log(\gamma_3/\gamma_*) \leq \log(2(r+1))$. Consider now the regime $1 \leq r < \sqrt{2}$. The number of mistakes for \cref{alg:max-margin} is at most $2+ \log(\gamma_3/\gamma_*)/(-\log(r/2)) \leq 2 + \log(2(r+1))/\log(2/r) \leq 2 + 2\log(2(\sqrt{2}+1))/\log(2) \leq 7$. Compared to the bound $O(r^2(\barD/D)^2)$, the main difference is again based on the size of $\barD/D$. This can be arbitrarily bad.

Finally, when $r \geq \sqrt{2}$, we have $0<1/r^2 <1$, and recall that the relation $z\le -\log(1-z)$ holds for all $0<z<1$, and thus $1/r^2 \leq -\log(1-1/r^2)$, implying $1/(-\log(1-1/r^2)) \le r^2$. Then, the mistake bound of \cref{alg:max-margin} is $2+ 2\log(\gamma_3/\gamma_*)/(-\log(1-1/r^2)) \le 2+ 2r^2 \log(2(r+1)) =  O(r^2 \log(2(r+1))$. Compared to $O(r^2(\barD/D)^2)$, in terms of $r$ this is worse by a logarithmic term, yet removes the translation-dependent term $\barD/D$. Therefore, if $\barD/D$ is large relative to $\log(2(r+1))$, then our bound is an improvement over the classical bound from the literature; and if $\barD/D$ is small compared to $\log(2(r+1))$, our bound is worse than the classical one by only a logarithmic term.

Note that it is well-known that the $O(r^2 (\barD/D)^2)$ mistake bound is tight for the perceptron, evidenced by a carefully constructed example \citep[Exercise 8.3]{mohri_foundations_2012}. Furthermore, it can be shown that \emph{any} deterministic algorithm also makes $d$ mistakes on the same example, where $d$ is the dimension of the problem, i.e., the feature vectors $x_t\in\bbR^d$. To force $d$ mistakes, the example constructs a sequence of $d$ feature vectors with $D = \sqrt{2}$, $\barD=1$ and $\gamma_* = 1/\sqrt{d}$. The mistake bound for \cref{alg:max-margin} in \cref{tab:perf_guarantees} does not contradict this: since $r = \sqrt{2d} \geq \sqrt{2}$, the mistake bound of \cref{alg:max-margin} is $4 \log(2) d \log(2(\sqrt{2d}+1)) \geq d$.

\section{Full experimental details}\label{sec:experiment-details}

In this section, we present full details of the experiments from \cref{sec:numerical}. Computations are performed on a personal laptop with Apple M3 processor and 24GB of memory, using Python and numpy. Conic optimization subproblems are solved with Mosek version 11.0.27.

\subsection{Preprocessing}\label{sec:exp-preprocessing}
Starting from the raw Adult dataset \citep{adult_dataset}, observations with missing features are removed and categorical variables with $k$ classes are converted into $k-1$ dummy variables. The dimension of the resulting dataset is $d = 96$. Each feature is also standardized by subtracting the respective sample mean and dividing by the standard deviation. Since the guarantees hold for separable data streams, we ensure separability by first fitting a linear support vector classifier on this dataset and then removing data points that are misclassified or within $0.01$ distance of the decision boundary. The resulting dataset consists of $35,508$ observations, with 27.11\% of the data having label +1.

\subsection{Data transformations}\label{sec:exp-transformations}
We perform three transformations of the dataset to test the effect of changing $\barD/D$ and $b_*$ on the performance of the algorithms. First, we compute the true maximum margin classifier $(w_*, b_*)$ on the preprocessed dataset described in \cref{sec:exp-preprocessing}, together with the certificate points $v_+,v_-$ as in \cref{lemma:separable-max-margin}. Let  $\gamma$ be the margin of this dataset, and let $\bar{x}$ be a feature vector in the dataset with the largest $\ell_2$-norm, i.e., $\barD = \max_{t}\|x_t\|_2 = \|\bar{x}\|_2$.

In our first transformation (margin normalization) we translate each point $(x,y) \to \left(x+\zeta \frac{y w_*}{\|w_*\|_2}, y\right)$ for fixed $\zeta \in \bbR$. This then changes the margin to $\gamma + \zeta$. We set $\zeta = 1-\gamma$, so that the margin of the new dataset is normalized to $\gamma_* = 1$, though the maximum margin classifier remains the same, i.e., $(w_*,b_*)$. We perform this transformation simply for convenience; $D$ is unaffected, and $r = D/\gamma_* = D$, though $\barD$ may change.

Our second transformation (bias control) affects the bias term $b_*$. In this transformation, we translate each point in the dataset $(x,y) \to \left(x - \frac{1}{2}(v_+ + v_-), y \right)$, where $v_+,v_-$ are the certificate points as in \cref{lemma:separable-max-margin}. After this transformation the new maximum margin classifier becomes $(w_*, 0)$, i.e., it has zero bias term. Note that parameters $r$ and $\gamma$ are unaffected by this transformation, though $\barD$ may change. This transformation is optional, and we test the effect with and without performing it.

Finally, our third transformation (feature scaling) scales the feature vectors based on an input parameter of $\theta$. In particular, we translate each point $(x,y) \to \left( x + \theta \left(\bar{x} - \frac{{\bar{x}}^\top w_*}{\|w_*\|_2^2} w_* \right), y \right)$ for a given $\theta \geq 0$ (recall $\barD = \|\bar{x}\|_2$). This transformation does not change $(w_*,b_*)$, $\gamma_*$ or $r$, but ensures that the new $\barD \geq \sqrt{\|\bar{x}\|_2^2 + \theta^2 \left\| \bar{x} - \frac{\bar{x}^\top w_*}{\|w_*\|_2^2} w_* \right\|_2^2}$. 

\subsection{Implementation details and hyper-parameter tuning}\label{sec:exp-setup}

We choose the tuning parameters $\phi, B$ and $C$ for \alma as described by \citet[Theorem 3]{gentile_new_2000}, i.e., we set $\phi=0.3$, $B = \sqrt{8}/(1-\phi)$, and $C=\sqrt{2}$. 
The \pumma algorithm requires a $\delta \in [0,1]$ parameter that has a similar interpretation to our aggressiveness parameter $\rho$: for $\delta \in (0,1)$ \citet{IshibashiEtAl2008} show that \pumma recovers $(1-\delta)$-fraction of the maximum margin, and setting $\delta$ near $0$ increases aggressiveness of \pumma. We choose $\delta = 0.01$ to enable \pumma to recover 99\% of the margin, and because \citet{IshibashiEtAl2008} have no guarantees for $\delta = 0$.

While \cref{alg:max-margin} (and thus \nomm and \eomm) and \pumma can handle non-zero bias $b_*$ automatically, the other algorithms require some adjustments to handle the case of $b_*\neq0$. \alma and the \pctr append a coordinate of $1$ to the feature vectors, thus increasing the dimension of the data by $+1$, with the bias term being the coefficient of $w$ corresponding to this last coordinate. On the other hand, this trick does not work for the \romma variants, so \citet{li_relaxed_1999} suggest instead adding a $-\barD$ coordinate, with the corresponding coefficient interpreted as approximating $-b_*/\barD$. This approach requires prior knowledge of $\barD$, and in our implementation we provided this information to both \romma and \aromma as input.

\subsection{Additional numerical results}\label{sec:extra_numerical_results}

We tested the impact of going through the same dataset multiple times. In particular, we allow all algorithms, except \nomm, to have five full passes through the data (each with a new random permutation of the points). In the case of \nomm only a single pass through the dataset is done, since \nomm will stop updating after seeing all the data once. We report our results in \cref{tab:perf}. In this second set of experiments, the margin guarantee of \eomm increases to 0.9330 while its number of mistakes remains stable at 6.
In the most favorable setting of $b_* = 0$ and $\theta = 0$, the margin performances of \aromma and \pumma improve slightly, yet is still lower than that of \eomm. When $b_* = 0$ and $\theta = 0$, the \pctr is able to achieve positive margin (meaning it has found a correct classifier) yet the margin is significantly lower than \eomm. Moreover, once again whenever $b_*\neq0$ or $\theta$ increases, the margin performances of \romma, \aromma and \alma degrade sharply and their numbers of mistakes increase significantly. On the other hand, the number of mistakes that the \pctr makes remains relatively stable, yet the margins achieved are significantly lower than \eomm, whenever they are positive.

\subsection{Translation invariance of \pumma}\label{sec:translation-invariance-pumma}
The \pumma algorithm has a similar initialization phase as \cref{alg:max-margin} in order to get one positive and one negative point, and the first $(w_1,b_1)$ is obtained by solving the maximum margin problem on those two points. This is translation invariant for a reason similar to \cref{rem:translation-invariance-offline}: only $b_1$ changes if a vector $u$ is added to both the positive point and the negative point. In each subsequent iteration, if an update occurs, it is computed by replacing either the positive point or the negative point with the new point $x_t$, according to the label $y_t$, then solving a maximum margin problem on the two points but with an additional constraint $w^\top w_{t-1} \geq \|w_{t-1}\|_2^2$. Since $w_{t-1}$ remains the same on translation, if a constant vector $u$ is added to the positive and negative points, the new $(w_t,b_t)$ only differs in the $b_t$ argument.

\begin{table}[htbp]
    \centering
    \resizebox{.95\textwidth}{!}{
    {\fontsize{10pt}{13.1pt}\selectfont
    \begin{tabular*}{\linewidth}{@{\extracolsep{\fill}}l|rrrr||rrrr}
    \toprule
        & \multicolumn{4}{c}{$b_* = 0$} & \multicolumn{4}{c}{$b_* = -0.2848$}\\
    \cline{2-5} \cline{6-9} 
        & $m$ & $\bar{\gamma}$ & $\tau$ & time (s) & $m$ & $\bar{\gamma}$ & $\tau$ & time (s) \\
    \midrule
    $\theta = 0$ & \multicolumn{4}{c}{$\barD/D = 0.9666$} & \multicolumn{4}{c}{$\barD/D = 0.9658$} \\
    \hline
    \nomm & 4 & 1.00 & 70 & 76.56 & 4 & 1.00 & 70 & 74.25 \\
    \ceomm & 21 & 0.05 & 16808 & 0.10 & 21 & 0.05 & 16808 & 0.11 \\
    \eomm & 6 & 0.84 & 534 & 0.11 & 6 & 0.84 & 534 & 0.12 \\
    \romma & 38 & 0.16 & 28065 & 0.09 & 4114 & -- & -- & 0.22 \\
    \aromma & 11 & 0.73 & 265 & 0.15 & 2643 & -- & -- & 0.21 \\
    \pumma & 9 & 0.68 & 584 & 0.08 & 9 & 0.68 & 584 & 0.08 \\
    \pctr & 59 & -- & -- & 0.07 & 52 & -- & -- & 0.07 \\
    \alma & 53 & -- & -- & 0.07 & 174 & -- & -- & 0.07 \\
    \midrule
    $\theta = 0.25$ & \multicolumn{4}{c}{$\barD/D = 1.2082$} & \multicolumn{4}{c}{$\barD/D = 1.2072$} \\
    \hline
    \nomm & 4 & 1.00 & 70 & 75.07 & 4 & 1.00 & 70 & 75.86 \\
    \ceomm & 21 & 0.05 & 16808 & 0.11 & 21 & 0.05 & 16808 & 0.11 \\
    \eomm & 6 & 0.84 & 534 & 0.12 & 6 & 0.84 & 534 & 0.12 \\
    \romma & 529 & -- & -- & 0.09 & 5600 & -- & -- & 0.30 \\
    \aromma & 220 & -- & 12282 & 0.16 & 3816 & -- & -- & 0.22 \\
    \pumma & 9 & 0.68 & 584 & 0.08 & 9 & 0.68 & 584 & 0.08 \\
    \pctr & 70 & -- & -- & 0.07 & 76 & -- & -- & 0.07 \\
    \alma & 689 & -- & 23099 & 0.08 & 1961 & -- & -- & 0.10 \\
    \midrule
    $\theta = 0.50$ & \multicolumn{4}{c}{$\barD/D = 1.4498$} & \multicolumn{4}{c}{$\barD/D = 1.4486$} \\
    \hline
    \nomm & 4 & 1.00 & 70 & 75.05 & 4 & 1.00 & 70 & 75.57 \\
    \ceomm & 21 & 0.05 & 16808 & 0.10 & 21 & 0.05 & 16808 & 0.11 \\
    \eomm & 6 & 0.84 & 534 & 0.12 & 6 & 0.84 & 534 & 0.13 \\
    \romma & 1571 & -- & -- & 0.12 & 6693 & -- & -- & 0.35 \\
    \aromma & 788 & -- & -- & 0.19 & 5096 & -- & -- & 0.21 \\
    \pumma & 9 & 0.68 & 584 & 0.08 & 9 & 0.68 & 584 & 0.09 \\
    \pctr & 66 & -- & -- & 0.07 & 56 & 0.03 & 24953 & 0.07 \\
    \alma & 2489 & -- & -- & 0.10 & 5448 & -- & -- & 0.12 \\
    \midrule
    $\theta = 0.75$ & \multicolumn{4}{c}{$\barD/D = 1.6914$} & \multicolumn{4}{c}{$\barD/D = 1.6901$} \\
    \hline
    \nomm & 4 & 1.00 & 70 & 74.63 & 4 & 1.00 & 70 & 75.74 \\
    \ceomm & 21 & 0.05 & 16808 & 0.11 & 21 & 0.05 & 16808 & 0.11 \\
    \eomm & 6 & 0.84 & 534 & 0.12 & 6 & 0.84 & 534 & 0.12 \\
    \romma & 2887 & -- & -- & 0.17 & 8099 & -- & -- & 0.43 \\
    \aromma & 1695 & -- & -- & 0.17 & 6232 & -- & -- & 0.22 \\
    \pumma & 9 & 0.68 & 584 & 0.08 & 9 & 0.68 & 584 & 0.08 \\
    \pctr & 112 & -- & -- & 0.07 & 74 & -- & -- & 0.07 \\
    \alma & 4633 & -- & -- & 0.11 & 8602 & -- & -- & 0.13 \\
    \midrule
    $\theta = 1$ & \multicolumn{4}{c}{$\barD/D = 1.9331$} & \multicolumn{4}{c}{$\barD/D = 1.9315$} \\
    \hline
    \nomm & 4 & 1.00 & 70 & 74.64 & 4 & 1.00 & 70 & 76.27 \\
    \ceomm & 21 & 0.05 & 16808 & 0.11 & 21 & 0.05 & 16808 & 0.11 \\
    \eomm & 6 & 0.84 & 534 & 0.12 & 6 & 0.84 & 534 & 0.12 \\
    \romma & 4132 & -- & -- & 0.22 & 9213 & -- & -- & 0.50 \\
    \aromma & 2728 & -- & -- & 0.18 & 7413 & -- & -- & 0.22 \\
    \pumma & 9 & 0.68 & 584 & 0.08 & 9 & 0.68 & 584 & 0.08 \\
    \pctr & 80 & -- & -- & 0.07 & 46 & 0.05 & 20579 & 0.07 \\
    \alma & 6636 & -- & -- & 0.11 & 10946 & -- & -- & 0.13 \\
    \bottomrule
    \end{tabular*}}
    }
    \caption{
    Numerical results over the dataset variants when all of the algorithms are allowed a single pass over the dataset.
    }
\end{table}

\begin{table}[htbp]
    \centering
    \resizebox{.95\textwidth}{!}{
    {\fontsize{10pt}{13.1pt}\selectfont
    \begin{tabular*}{\linewidth}{@{\extracolsep{\fill}}l|rrrr||rrrr}
    \toprule
        & \multicolumn{4}{c}{$b_* = 0$} & \multicolumn{4}{c}{$b_* = -0.2848$}\\
    \cline{2-5} \cline{6-9} 
        & $m$ & $\bar{\gamma}$ & $\tau$ & time (s) & $m$ & $\bar{\gamma}$ & $\tau$ & time (s) \\
    \midrule
    $\theta = 0$ & \multicolumn{4}{c}{$\barD/D = 0.9666$} & \multicolumn{4}{c}{$\barD/D = 0.9658$} \\
    \hline
    \nomm & 4 & 1.00 & 70 & 76.56 & 4 & 1.00 & 70 & 74.25 \\
    \ceomm & 21 & 0.05 & 16808 & 0.20 & 21 & 0.05 & 16808 & 0.22 \\
    \eomm & 6 & 0.93 & 534 & 0.55 & 6 & 0.93 & 534 & 0.59 \\
    \romma & 38 & 0.16 & 28065 & 0.18 & 6637 & -- & -- & 1.17 \\
    \aromma & 11 & 0.89 & 265 & 0.80 & 3095 & -- & 150864 & 1.06 \\
    \pumma & 9 & 0.88 & 584 & 0.37 & 9 & 0.88 & 584 & 0.40 \\
    \pctr & 69 & 0.06 & 64028 & 0.21 & 78 & 0.01 & 170119 & 0.38 \\
    \alma & 79 & -- & -- & 0.35 & 329 & -- & -- & 0.35 \\
    \midrule
    $\theta = 0.25$ & \multicolumn{4}{c}{$\barD/D = 1.2082$} & \multicolumn{4}{c}{$\barD/D = 1.2072$} \\
    \hline
    \nomm & 4 & 1.00 & 70 & 75.07 & 4 & 1.00 & 70 & 75.86 \\
    \ceomm & 21 & 0.05 & 16808 & 0.22 & 21 & 0.05 & 16808 & 0.22 \\
    \eomm & 6 & 0.93 & 534 & 0.58 & 6 & 0.93 & 534 & 0.60 \\
    \romma & 577 & -- & -- & 0.50 & 9531 & -- & -- & 1.55 \\
    \aromma & 221 & 0.56 & 12282 & 0.85 & 4851 & -- & -- & 1.10 \\
    \pumma & 9 & 0.88 & 584 & 0.39 & 9 & 0.88 & 584 & 0.41 \\
    \pctr & 72 & 0.16 & 60814 & 0.22 & 122 & -- & -- & 0.39 \\
    \alma & 695 & 0.64 & 23099 & 0.36 & 2142 & 0.36 & 45381 & 0.41 \\
    \midrule
    $\theta = 0.50$ & \multicolumn{4}{c}{$\barD/D = 1.4498$} & \multicolumn{4}{c}{$\barD/D = 1.4486$} \\
    \hline
    \nomm & 4 & 1.00 & 70 & 75.05 & 4 & 1.00 & 70 & 75.57 \\
    \ceomm & 21 & 0.05 & 16808 & 0.22 & 21 & 0.05 & 16808 & 0.23 \\
    \eomm & 6 & 0.93 & 534 & 0.57 & 6 & 0.93 & 534 & 0.62 \\
    \romma & 2039 & -- & -- & 0.62 & 13052 & -- & -- & 2.15 \\
    \aromma & 807 & 0.07 & 39063 & 0.84 & 7292 & -- & -- & 1.09 \\
    \pumma & 9 & 0.88 & 584 & 0.38 & 9 & 0.88 & 584 & 0.41 \\
    \pctr & 70 & 0.00 & 47354 & 0.22 & 56 & 0.03 & 24953 & 0.15 \\
    \alma & 2685 & 0.28 & 57793 & 0.39 & 7999 & 0.01 & 135655 & 0.47 \\
    \midrule
    $\theta = 0.75$ & \multicolumn{4}{c}{$\barD/D = 1.6914$} & \multicolumn{4}{c}{$\barD/D = 1.6901$} \\
    \hline
    \nomm & 4 & 1.00 & 70 & 74.63 & 4 & 1.00 & 70 & 75.74 \\
    \ceomm & 21 & 0.05 & 16808 & 0.22 & 21 & 0.05 & 16808 & 0.23 \\
    \eomm & 6 & 0.93 & 534 & 0.58 & 6 & 0.93 & 534 & 0.62 \\
    \romma & 4140 & -- & -- & 0.83 & 16850 & -- & -- & 2.85 \\
    \aromma & 1829 & -- & 89998 & 0.91 & 10124 & -- & -- & 1.09 \\
    \pumma & 9 & 0.88 & 584 & 0.39 & 9 & 0.88 & 584 & 0.40 \\
    \pctr & 138 & -- & -- & 0.37 & 82 & 0.01 & 52263 & 0.23 \\
    \alma & 5926 & -- & -- & 0.43 & 14866 & -- & -- & 0.49 \\
    \midrule
    $\theta = 1$ & \multicolumn{4}{c}{$\barD/D = 1.9331$} & \multicolumn{4}{c}{$\barD/D = 1.9315$} \\
    \hline
    \nomm & 4 & 1.00 & 70 & 74.64 & 4 & 1.00 & 70 & 76.27 \\
    \ceomm & 21 & 0.05 & 16808 & 0.22 & 21 & 0.05 & 16808 & 0.22 \\
    \eomm & 6 & 0.93 & 534 & 0.59 & 6 & 0.93 & 534 & 0.59 \\
    \romma & 6653 & -- & -- & 1.19 & 20737 & -- & -- & 3.58 \\
    \aromma & 3176 & -- & 160346 & 0.92 & 13393 & -- & -- & 1.13 \\
    \pumma & 9 & 0.88 & 584 & 0.40 & 9 & 0.88 & 584 & 0.39 \\
    \pctr & 84 & 0.09 & 45807 & 0.22 & 46 & 0.05 & 20579 & 0.15 \\
    \alma & 9812 & -- & -- & 0.48 & 23840 & -- & -- & 0.53 \\
    \bottomrule
    \end{tabular*}}
    }
    \caption{
    Numerical results over the dataset variants when all of the algorithms (except \nomm) are allowed five passes over the dataset.
    }
    \label{tab:perf}
\end{table}

\end{document}